\documentclass[11pt]{article}

\usepackage{amsfonts}

\usepackage{amsmath}

\newcommand{\beq}{\begin{equation}}

\newcommand{\eeq}{\end{equation}}

\newcommand{\bea}{\begin{eqnarray}}

\newcommand{\eea}{\end{eqnarray}}

\newcommand{\beas}{\begin{eqnarray*}}

\newcommand{\eeas}{\end{eqnarray*}}

\newtheorem{theorem}{Theorem}[section]

\newtheorem{definition}[theorem]{Definition}

\newtheorem{proposition}[theorem]{Proposition}

\newtheorem{lemma}[theorem]{Lemma}

\newtheorem{remark}[theorem]{Remark}

\newtheorem{example}[theorem]{Example}

\newtheorem{examples}[theorem]{Examples}

\newtheorem{foo}[theorem]{Remarks}

\newenvironment{proof}{\addvspace{\medskipamount}\par\noindent{\it

Proof}.}{\unskip\nobreak\hfill$\Box$\par\addvspace{\medskipamount}}

\title{Self-similarity and fractional Brownian motions on Lie groups}

\author{Fabrice Baudoin
\\\small Laboratoire de Statistiques et Probabilit\'es
\\\small Universit\'e Paul Sabatier
\\\small 118 Route de Narbonne, Toulouse, France
\\\small fbaudoin@cict.fr
\and Laure Coutin
\\\small Laboratoire MAP5
\\\small Universit\'e Ren\'e Descartes
\\\small 45, rue des Saints P\`eres, Paris, France
\\\small coutin@cict.fr }

\begin{document}

\maketitle

\begin{abstract}

The goal of this paper is to define and study a notion of
fractional Brownian motion on a Lie group. We define it as at the
solution of a stochastic differential equation driven by a linear
fractional Brownian motion. We show that this process has
stationary increments and satisfies a local self-similar property.
Furthermore the Lie groups for which this self-similar property is
global are characterized. Finally, we prove an integration by
parts formula on the path group space and deduce the existence of
a density.

\

\

\

\end{abstract}

\tableofcontents

\newpage

\section{Introduction}
Since the seminal works of It\^o \cite{Ito}, Hunt \cite{Hunt}, and
Yosida \cite{Yosida}, it is well known that the (left) Brownian
motion on a Lie group $\mathbf{G}$ appears as the solution of a
stochastic differential equation
\begin{align}\label{SDEliegroup}
dX_t = \sum_{i=1}^d  V_i (X_t) \circ dB^i_t, \text{ }t \geq 0,
\end{align}
where $V_1,...,V_d$ are left-invariant vector fields on
$\mathbf{G}$ and where $(B_t)_{t \ge 0}$ is a Brownian motion; for
further details on this, we also refer to \cite{Bau} and
\cite{Ro-Wi}. In this paper we investigate the properties of the
solution of an equation of the type (\ref{SDEliegroup}) when the
driving Brownian motion is replaced by a fractional Brownian
motion with parameter $H$. Recently, there have been numerous
attempts to define a notion of solution for differential equations
driven by fractional Brownian motion. One dimensional differential
equation can be solved  using a Doss-Susmann aproach for any
values of the parameter $H$ as in the work of Nourdin,
\cite{nourdin} or in the linear case by Nourdin-Tudor in
\cite{nourdin-tudor} .  The situation is quite different in the
multidimensional case. When the Hurst parameter is greater than
$1/2$ existence and uniqueness of the solution are obtained by
Z\"ahle in \cite{zahle} or Nualart-Rasc\^anu in \cite{rascanu}.
And finally, as a consequence of the work of Coutin and Qian
\cite{CQ} a notion of solution is actually well-defined for $H
>\frac{1}{4}$.

\

The paper is organized as follows.

In a first section, we show existence and uniqueness for the
solution of an equation of the type (\ref{SDEliegroup}) when
$(B_t)_{t \ge 0}$ is a fractional Brownian motion with parameter
$H >\frac{1}{4}$. The solution is shown to have stationary
increments. We also check that the solution in invariant in law by
isometries.

In the second section, we study the scaling properties of the
solution. In the spirit of the notion of asymptotic
self-similarity studied by Kunita \cite{Kunita1}, \cite{Kunita2},
(see also \cite{cohen}), we show that the fractional Brownian
motion on the group is asymptotically self-similar with parameter
$H$. After that, we characterize the groups for which the scaling
property is global: such groups are necessarily simply connected
and nilpotent.

Finally, the goal of the last section is to show the existence of
a density for the solution. At this end, we prove an integration
by parts formula on the path group space.

\

To simplify the presentation of our results, we mainly worked in
the setting of Lie groups of matrices. Nevertheless all our
results extend to general Lie groups.

\newpage

\section{Fractional Brownian motion on a Lie group}

Let us first recall that a $d$-dimensional fractional Brownian
motion with Hurst parameter $H \in (0,1)$ is a Gaussian process
\[
B_t = (B_t^1,...,B_t^d), \text{ } t \geq 0,
\]
where $B^1,...,B^d$ are $d$ independent centered Gaussian
processes with covariance function
\[
R\left( t,s\right) =\frac{1}{2}\left(
s^{2H}+t^{2H}-|t-s|^{2H}\right).
\]
It can be shown that such a process admits a continuous version
whose paths are H\"older $p$ continuous, $p<H$. Let us observe
that for $H= \frac{1}{2}$, $B$ is a Brownian motion.

Let $\mathbf{G}$ be a finite-dimensional ($\dim \mathbf{G} =d$)
connected Lie group of matrices with Lie algebra $\mathfrak{g}$.
We consider a basis $(V_1,...,V_d)$ of $\mathfrak{g}$.

If $(B_t^1,...,B_t^d)_{t \geq 0}$ is a $d$-dimensional fractional
Brownian motion in $\mathbb{R}^d$ with Hurst parameter $H \in
(0,1)$. The process
\[
B^{\mathfrak{g}}_t=\sum_{i=1}^d B^i_t V_i
\]
shall be called the canonical fractional Brownian motion on
$\mathfrak{g}$ with respect to the basis $(V_1,...,V_d)$.

In the remainder of this section, we assume now $H>\frac{1}{4}$.





\begin{theorem} \label{fbmliegroup} The equation
\begin{align}\label{eq:ed-def}
d X_t =X_t dB^{\mathfrak{g}}_t,~~~~ X_0=1_{\mathbf{G}}
\end{align}
has a unique solution in $\mathbf{G}$ in the sense of rough paths
of \cite{LyQi}. This solution $(X_t)_{t \ge 0}$ satisfies for
every $ s \ge 0$, $(X_s^{-1}X_{t+s} )_{t \ge 0}=^{law}(X_t)_{t \ge
0}$. The process $(X_t)_{t \ge 0}$ shall be called a  left
fractional Brownian motion with parameter $H$ on $\mathbf{G}$ with
respect to the basis $(V_1,...,V_d).$
\end{theorem}

\begin{proof}
We first show the existence and the uniqueness of a solution in
$\mathbf{G}$. Without loss of generality, we can work on the time
interval $[0,1]$. First, according to \cite{CQ}, the equation
(\ref{eq:ed-def}) has a unique solution with finite $p$ variation
for $p >\frac{1}{H}.$ Secondly, let $B^{\mathfrak{g},m}$ be the
sequel of linear interpolation of $B^{\mathfrak{g}}$ along the
dyadic subdivision of mesh $m$; that is if $t_i^m= i 2^{-m}$ for
$i=0,..., 2^m;$ then for $t \in [t_i^m, t_{i+1}^m ),$
\begin{align*}
B^{\mathfrak{g},m}_t=B^{\mathfrak{g}}(t_{i^m}) +
\frac{t-t_{i^m}}{t_{i+1}^m - t_i^m} (B^{\mathfrak{g}}_{t_{i+1}^m}
-B^{\mathfrak{g}}_{t_{i}^m} ).
\end{align*}
Let us now denote $X^m$ the solution of (\ref{eq:ed-def}) where
$B^{\mathfrak{g}}$ is replaced by $B^{\mathfrak{g},m}$, that is

\[
dX^m_t =\sum_{k=1}^d X^m_t      dB^{k,m}_tV_k
\]
It is easily seen that for $t \in [t_{n-1}^m, t_{n}^m )$,
recursevely on $n,$ $n=,...,2^m-1,$
\[
X^m_t =\exp \left(2^m (t-t^m_{n-1})\sum_{k=1}^d (B^k_{t^m_n}
-B_{t^m_{n-1}}^k )V_k \right) \cdots \exp \left( \sum_{k=1}^d
(B^k_{t^m_1}-B_{t^m_{0}}^k )V_k \right).
\]
Therefore, $X^m$ takes its values in $\mathbf{G}$. Now, from
\cite{CQ},  $X^m$ converges to $X$ for the distance of $1/p$
H\"older. Since the group $\mathbf{G}$ is closed in the $1/p$
H\"older topology, we conclude that $X$ belongs to  $\mathbf{G}.$

We now show that for every $s \ge 0$, the processes
$((X_s)^{-1}X_{t+s} ,~~t \geq 0)$ and $(X_t  ,~~t \geq 0)$ have
the same law. Let us fix $s \ge 0$. Once time again, the idea is
to use a linear interpolation along the dyadic subdivision of
$[0,1]$ of mesh $m$ and we keep the previous notations. First, let
us observe that for $t \ge s$,
\begin{align}\label{liegroupst}
X_{t} =X_s+\int_s^{t} X_u dB^{\mathfrak{g}}_u.
\end{align}
Let us now denote $(X^{m,s}_t)_{s \le t \le s+1}$ the solution of
(\ref{liegroupst}) where $B^{\mathfrak{g}}$ is replaced by
$(B^{\mathfrak{g},m}_{s+t})_{0 \le t \le 1}$. Therefore, for $t
\in [t_{n-1}^m, t_{n}^m )$,
\[
X^{m,s}_{t+s} =X_s \exp \left(2^m (t-t_{n-1}^m)\sum_{k=1}^d
(B^k_{s+t^m_{n}} -B_{s+t^m_{n-1}}^k )V_k \right) \cdots \exp
\left(\sum_{k=1}^d (B^k_{s+t^m_1} -B_{s+t^m_{0}}^k )V_k \right).
\]
By using the stationarity of the increments of the Euclidean
fractional Brownian motion, we get therefore:
\[
(X_s^{-1}X^{m,s}_{t+s} )_{0 \le t \le 1}=^{law} (X^m_t)_{0 \le t
\le 1}
\]
Using the Wong-Zakai theorem of \cite{CQ} and passing to the
limit, we obtain that for every $ s \ge 0$, $(X_s^{-1}X_{t+s}
)_{t\ge 0}=^{law}(X_t)_{t \ge 0}$.
\end{proof}

\begin{remark}
In the same way, we call the solution of the differential equation
\[
d X_t =dB^{\mathfrak{g}}_t X_t ,~~~~ X_0=\mathbf{1}_{\mathbf{G}}.
\]
a right fractional Brownian motion on $\mathbf{G}$. It is easily
seen that if $(X_t)_{t \ge 0}$ is a left fractional Brownian
motion on $\mathbf{G}$, then $(X^{-1}_t)_{t \ge 0}$ is a right
fractional Brownian motion on $\mathbf{G}$.
\end{remark}

Let us now turn to some  examples.

\begin{example}
The first basic example is $(\mathbf{R}^d,+)$. In that case, the
Lie algebra is generated by the vector fields
$\frac{\partial}{\partial x_1}, ...,\frac{\partial}{\partial x_d}$
and the fractional Brownian on $(\mathbf{R}^d,+)$ is nothing else
but the usual Euclidean fractional Brownian motion.
\end{example}

\begin{example}
The second basic example is the circle. Let
\[
\mathbf{S}^1=\left\{ z \in \mathbb{C}, \mid z \mid =1 \right\}.
\]
The Lie algebra of $\mathbf{S}^1$ is $\mathbb{R}$ and is generated
by $\frac{\partial}{\partial \theta}$ and the fractional Brownian
motion on $\mathbf{S}^1$ is given by
\[
X_t=e^{i B_t}, \quad t \ge 0,
\]
where $(B_t)_{t \ge 0}$ is a fractional Brownian motion on
$\mathbb{R}$.
\end{example}

\begin{example}
Let us consider the Lie group $\mathbf{SO} (3)$, i.e. the group of
$3 \times 3$, real, orthogonal matrices of determinant $1$. Its
Lie algebra $\mathfrak{so} (3)$ consists of $3 \times 3$, real,
skew-adjoint matrices of trace $0$. A basis of $\mathfrak{so} (3)$
is formed by
\[
V_{1}=\left(
\begin{array}{ccc}
~0~ & ~1~ & ~0~ \\
-1~ & ~0~ & ~0~ \\
~0~ & ~0~ & ~0~
\end{array}
\right) ,\text{ }V_{2}=\left(
\begin{array}{ccc}
~0~ & ~0~ & ~0~ \\
~0~ & ~0~ & ~1~ \\
~0~ & -1~ & ~0~
\end{array}
\right) ,\text{ }V_{3}=\left(
\begin{array}{ccc}
~0~ & ~0~ & ~1~ \\
~0~ & ~0~ & ~0~ \\
-1~ & ~0~ & ~0~
\end{array}
\right)
\]
A left fractional Brownian motion on $\mathbf{SO} (3)$ is
therefore given by the solution of the linear equation
\[
dX_t=X_t \left(
\begin{array}{ccc}
~0~ & ~dB^1_t~ & ~dB^3_t~ \\
~-dB^1_t~ & ~0~ & ~dB^2_t~ \\
~-dB^3_t~ & -dB^2_t~ & ~0~
\end{array}
\right), \quad X_0=1.
\]
\end{example}
This notion of fractional Brownian motion on a Lie group is
invariant by isometries, so that the law is invariant by an
orthonormal change of basis. More precisely, let us consider the
scalar product on $\mathfrak{g}$ that makes the basis
$V_1,...,V_d$ orthonormal. This scalar product defines a
Riemannian structure on $\mathbf{G}$ for which the left action is
an action by isometries. We have the following proposition:

\begin{proposition}
Let $\Psi: \mathbf{G} \rightarrow \mathbf{G}$ be a Lie group
morphism such that $d\Psi_{1_{\mathbf{G}}}$ (differential of
$\Psi$ at $1_{\mathbf{G}}$) is an isometry and let $(X_t)_{t \ge
0}$ be the left fractional Brownian motion on $\mathbf{G}$ as
defined in Theorem \ref{fbmliegroup}. We have:
\[
(\Psi(X_t))_{t \ge 0}=^{law} (X_t)_{t \ge 0}.
\]
\end{proposition}
\begin{proof}
Let us observe that by the change of variable formula, see
\cite{lyonsqian}:
\[
d \Psi (X_t) =\Psi(X_t) \left( \sum_{i=1}^d
d\Psi_{1_{\mathbf{G}}}(V_i) B^i_t \right),
\]
Now,
\[
\left( \sum_{i=1}^d d\Psi_{1_{\mathbf{G}}}(V_i) B^i_t \right)_{t
\ge 0} =^{law}\left( \sum_{i=1}^d V_i B^i_t \right)_{t \ge 0},
\]
because of the orthogonal invariance of the Euclidean fractional
Brownian motion. Therefore,
\[
(\Psi(X_t))_{t \ge 0}=^{law} (X_t)_{t \ge 0}.
\]
\end{proof}

\begin{remark}
If $\mathbf{G}$ is compact then there exists a bi-invariant
Riemannian metric and so, if $(X_t)_{t \ge 0}$ denotes a left
fractional Brownian motion for this bi-invariant metric, from the
previous proposition, we get  that for every $g \in \mathbf{G}$,
\[
(gX_tg^{-1})_{t \ge 0}=^{law} (X_t)_{t \ge 0}.
\]
\end{remark}

If the group $\mathbf{G}$ is nilpotent then we have a closed
formula for the left fractional Brownian motion on $\mathbf{G}$
that extends the well-known formula for the Brownian motion on a
nilpotent group (see by e.g. \cite{Bau}, \cite{Cast} or
\cite{Yam}).

Let us introduce some notations: For $k \ge 1$,

\begin{itemize}
\item
\[
\Delta^k [0,t]=\{ (t_1,...,t_k) \in [0,t]^k, t_1 < ... < t_k \};
\]
\item If $I=(i_1,...i_k) \in \{1,...,d\}^k$ is a word with length
$k$,
\begin{equation*}
\int_{\Delta^k [0,t]}  dB^I= \int_{0 < t_1 < ... < t_k \leq t}
dB^{i_1}_{t_1}  ...  dB^{i_k}_{t_k};
\end{equation*}
\item We denote $\mathfrak{S}_k$ the group of the permutations of the
index set $\{1,...,k\}$ and if $\sigma \in \mathfrak{S}_k$, we
denote for a word $I=(i_1,...,i_k)$,  $\sigma \cdot I$ the word
$(i_{\sigma(1)},...,i_{\sigma(k)})$;
\item If
$I=(i_1,...,i_k) \in \{ 1,..., d \}^k$ is a word, we denote by
$V_I$ the Lie commutator defined by
\[
V_I = [V_{i_1},[V_{i_2},...,[V_{i_{k-1}}, V_{i_{k}}]...];
\]
\item If $\sigma \in \mathfrak{S}_k$, we denote $e(\sigma)$ the
cardinality of the set
\[
\{ j \in \{1,...,k-1 \} , \sigma (j) > \sigma(j+1) \};
\]
\item Finally, if $I=(i_1,...,i_k) \in \{ 1,..., d \}^k$ is a word
\[
\Lambda_I (B)_t=\sum_{\sigma \in \mathfrak{S}_k} \frac{\left(
-1\right) ^{e(\sigma )}}{k^{2}\left(
\begin{array}{l}
k-1 \\
e(\sigma )
\end{array}
\right) } \int_{\Delta^k [0,t]} \circ dB^{\sigma^{-1} \cdot I}.
\]
\end{itemize}

\begin{proposition}\label{fbmnilpotent}
Assume that $\mathbf{G}$ is a nilpotent group then:
\[
X_t=\exp \left( \sum_{k = 1}^{+\infty} \sum_{I=(i_1,...,i_k)}
\Lambda_I (B)_t V_I \right), \quad t \ge 0,
\]
where the above sum is actually finite and where $(X_t)_{t \ge 0}$
is the left fractional Brownian motion defined as in Theorem
\ref{fbmliegroup}.
\end{proposition}

\begin{proof}
Let $B^{m}$ be the sequel of linear interpolation of
$B^{\mathfrak{g}}$ along the dyadic subdivision of mesh $m$. Let
us now denote $X^m$ the solution of (\ref{eq:ed-def}) where $B$ is
replaced by $B^{m}$. As already seen, for $t \in [t_{n-1}^m,
t_{n}^m )$,
\[
X^m_t =\exp \left(2^m(t-t_{n-1}^m) \sum_{k=1}^d (B^k_t
-B_{t^m_{n-1}}^k )V_k \right) \cdots \exp \left( \sum_{k=1}^d
(B^k_{t^m_1} -B_{t^m_{0}}^k )V_k \right).
\]
Now we use the Baker-Campbell-Hausdorff formula in nilpotent Lie
groups (see \cite{Bau}, \cite{Cast}, \cite{Stri}) to write the
previous product of exponentials under the form
\[
X^m_t=\exp \left( \sum_{k=1}^{+\infty} \sum_{I=\{i_1,...,i_k\}}
\Lambda_I (B^m)_t V_I \right),
\]
where
\[
\Lambda_I (B^m)_t=\sum_{\sigma \in \mathfrak{S}_k} \frac{\left(
-1\right) ^{e(\sigma )}}{k^{2}\left(
\begin{array}{l}
k-1 \\
e(\sigma )
\end{array}
\right) } \int_{\Delta^k [0,t]} \circ dB^{m,\sigma^{-1} \cdot I}.
\]
From \cite{CQ}, in the distance of $p$ variation, with $p >
\frac{1}{H}$, and if the length of the word $I$ is less than  2,
\begin{align}\label{cv-int-2}
\Lambda_I (B^m)_t \rightarrow_{m \rightarrow +\infty} \Lambda_I
(B)_t.
\end{align}
By using now Theorem 3.1.3  of \cite{lyonsqian}, the convergence
in (\ref{cv-int-2}) holds for all word. Therefore,
\[
X_{t}= \exp \left( \sum_{k = 1}^{+\infty} \sum_{I=(i_1,...,i_k)}
\Lambda_I (B)_t V_I \right), \text{ } t \geq 0.
\]
\end{proof}

\begin{example}
In a two-step nilpotent group, that if is all brackets with length
more than two are zero, we have therefore
\[
X_t=\exp \left( \sum_{i=1}^d B^i_t V_i + \frac{1}{2} \sum_{1 \le i
<j \le d} \left( \int_0^t B^i_s dB^j_s -B^j_s dB^i_s \right)
[V_i,V_j] \right).
\]
\end{example}

\section{Self-similarity of a fractional Brownian motion on a Lie group}

Recall that for the Euclidean fractional Brownian motion, we have
\[
(B_{ct}^1,...,B_{ct}^d)_{t \ge 0}=^{law} (c^H B_{t}^1,...,c^H
B_{t}^d)_{t \ge 0}
\]
This property is called the scaling property of the fractional
Brownian motion. In this section, we are going to study scaling
properties of fractional Brownian motions on a Lie group.

As in the previous section, let $\mathbf{G}$ be a connected Lie
group (of matrices) with Lie algebra $\mathfrak{g}$. Let
$V_1,...,V_d$ be a basis of $\mathfrak{g}$ and denote by $(X_t)_{t
\ge 0}$ the solution of the equation
\[
dX_t= X_t \left( \sum_{i=1}^d V_i dB^i_t \right), \quad
X_0=1_{\mathbf{G}},
\]
where $(B_t)_{t \ge 0}$ is a $d$-dimensional fractional Brownian
with Hurst parameter $H > \frac{1}{3}$. We restrict to the case
$H>\frac{1}{3}$ to use some technical estimates that come from
\cite{Baudoin-Coutin} but the results of this section certainly
also hold for $H>\frac{1}{4}$.

\subsection{Local self-similarity}

First, we notice that for $(X_t)_{t \ge 0}$ we always have an
asymptotic scaling property in the following sense:
\begin{proposition}\label{localselfsimilarity}
Let $f:\mathbf{G} \rightarrow \mathbb{R}$ be a smooth map such
that $ \sum_{i=1}^d (V_i f)(1_{\mathbf{G}})^2 \neq 0$. Then, when
$c \rightarrow 0$, $c>0$, the sequence of processes
$\left(\frac{1}{c^H}\left( f(X_{ct})-f(1_{\mathbf{G}}) \right)
\right)_{0 \le t \le 1}$ converges in law to   $(a \beta_t)_{0 \le
t \le 1}$ where $(\beta_t)_{t \ge 0}$ is a one-dimensional
fractional Brownian motion and
\[
a=\sqrt{ \sum_{i=1}^d (V_i f)(1_{\mathbf{G}})^2}.
\]
\end{proposition}

\begin{proof}
From   \cite{Baudoin-Coutin} inequality (4.7),
\[
f(X_t)= f(1_{\mathbf{G}}) + \sum_{i=1}^d (V_i f)
(1_{\mathbf{G}})B^i_t+\mathbf{R} (t), \text{ } t \geq 0
\]
for some remainder term $\mathbf{R}$ that satisfies
\begin{align*}
|\mathbf{R} (t)| \leq Ct^{2/p}
\end{align*}
where $C$ is a random variable with finite exponential moment and
$p > 1/H$. Therefore,
\[
\frac{1}{c^H}\left( f(X_{ct})-f(1_{\mathbf{G}})
\right)=\sum_{i=1}^d (V_i f)
(1_{\mathbf{G}})\frac{B^i_{ct}}{c^H}+\frac{\mathbf{R}
(ct)}{c^H},\text{ } t\ge 0,
\]
and the convergence follows easily from the scaling property of
the fractional Brownian motion.
\end{proof}

\begin{remark}
Slightly more generally, by using the Taylor expansion proved in
\cite{Baudoin-Coutin}, we obtain in the same way: Let
$f:\mathbf{G} \rightarrow \mathbb{R}$ be a smooth map such that
there exist $k \ge 1$ and $(i_1,...,i_k) \in \{ 1,...,d \}^k$ that
satisfy
\[
(V_{i_1} \cdots V_{i_k} f)(1_{\mathbf{G}}) \neq 0.
\]
Denote $n$ the smallest $k$ that satisfies the above property.
Then, when $c \rightarrow 0$, $c>0$, the sequence of processes
$\left(\frac{1}{c^{nH}}\left( f(X_{ct})-f(1_{\mathbf{G}}) \right)
\right)_{0 \le t \le 1}$ converges in law to   $( \beta_t)_{0 \le
t \le 1}$ where $(\beta_t)_{t \ge 0}$ is such that
\[
(\beta_{ct})_{t \ge 0} =^{law}(c^{nH}\beta_{t})_{t \ge 0}.
\]
\end{remark}

\subsection{Global self-similarity}

Despite the local self-similar property, as we will see, in
general there is no global scaling property for the fractional
Brownian motion on a Lie group. Let us first briefly discuss what
should be a good notion of scaling in a Lie group (see also
\cite{Kunita1} and \cite{Kunita2}). If we can find a family a map
$\Delta_c$, $c>0$, such that
\[
\left( X_{ct}\right)_{t \ge 0}=^{law}\left( \Delta_c
X_{t}\right)_{t \ge 0},
\]
first of all it is natural to require that the map $c \rightarrow
\Delta_c$ is continuous and $\lim_{c \rightarrow 0}
\Delta_c=1_{\mathbf{G}}$. Then by looking at $\left( X_{c_1 c_2
t}\right)_{t \ge 0}$, we will also naturally ask that $\Delta_{c_1
c_2} = \Delta_{c_1 } \circ \Delta_{c_2}$. Finally, since
$(X_s^{-1}X_{t+s} )_{t \ge 0}=^{law}(X_t)_{t \ge 0}$, we will also
ask that $\Delta_c$ is a Lie group automorphism.

The following theorem shows that the existence of such a family
$\Delta_c$ on $\mathbf{G}$ only holds if the group is
$(\mathbf{R}^d,+)$. This is partly due to the following lemma of
Lie group theory that says that the existence of a dilation on
$\mathbf{G}$ imposes strong topological and algebraic
restrictions:

\begin{lemma}\label{techniclemma}(See \cite{Lawton})
Assume that there exists a Lie group automorphism $\Psi :
\mathbf{G} \rightarrow \mathbf{G}$  such that
$d\Psi_{1_{\mathbf{G}}}$ (differential of $\Psi$ at
$1_{\mathbf{G}}$) has all its eigenvalues of modulus $>1$, then
$\mathbf{G}$ is a simply connected nilpotent Lie group.
\end{lemma}

We can now show:
\begin{theorem}\label{theoremescaling}
There exists a family of Lie group automorphisms $\Delta_t
:\mathbf{G} \rightarrow \mathbf{G}$, $t > 0$, such that:
\begin{enumerate}
\item The map $t \rightarrow \Delta_t$ is continuous and $\lim_{t \rightarrow 0} \Delta_t=1_{\mathbf{G}}$;
\item For $t_1,t_2 \ge 0$, $\Delta_{t_1 t_2} = \Delta_{t_1 } \circ \Delta_{t_2}$;
\item $(X_{ct})_{t \geq 0} =^{law} (\Delta_c X_{t})_{t \geq 0}$;
\end{enumerate}
if and only if the group $\mathbf{G}$ is isomorphic to
$(\mathbf{R}^d,+)$.
\end{theorem}

\begin{proof}

If $\mathbf{G}$ is isomorphic to $(\mathbf{R}^d,+)$, $(X_t)_{t \ge
0}$ is a Euclidean fractional Brownian motion and the result is
trivial.

We prove now the converse statement. Let us first show that the
existence of the family $(\Delta_t)_{t> 0}$ implies that
$\mathbf{G}$ is a simply connected nilpotent Lie group.

Let us denote
\[
\delta_c =d\Delta_c (1_{\mathbf{G}}),
\]
the differential map of $\Delta_c$ at $1_{\mathbf{G}}$ and observe
that $\delta_c$ is a Lie algebra automorphism $\mathfrak{g}
\rightarrow \mathfrak{g}$. The map $f: t \rightarrow \delta_{e^t}$
is a map from $\mathbf{R}$ onto the set of linear maps
$\mathfrak{g} \rightarrow \mathfrak{g}$ that is continuous. We
have furthermore the property
\[
f (t+s)= f (t) f (s).
\]
Consequently there exists a linear map $\phi : \mathfrak{g}
\rightarrow \mathfrak{g}$ such that
\[
\delta_c = \text{Exp} ( \phi \ln c), \text{   } c>0,
\]
where $\text{Exp}$ denotes here the exponential of linear maps
(and not the exponential map $\mathfrak{g} \rightarrow \mathbf{G}$
which is denoted $\exp$). Let us furthermore observe that if
$\lambda \in \mathbf{Sp} (\phi)$ is an eigenvalue of $\phi$, then
$\Re \phi >0$ because $\lim_{c \rightarrow 0}
\Delta_c=1_{\mathbf{G}}$. Therefore from Lemma \ref{techniclemma},
$\mathbf{G}$ has to be a simply connected nilpotent Lie group.

We deduce from Proposition \ref{fbmnilpotent} that
\[
X_t=\exp \left( \sum_{k = 1}^{+\infty} \sum_{I=(i_1,...,i_k)}
\Lambda_I (B)_t V_I \right), \quad t \ge 0,
\]
where the above sum is actually finite and
\[
\Lambda_I (B)_t=\sum_{\sigma \in \mathfrak{S}_k} \frac{\left(
-1\right) ^{e(\sigma )}}{k^{2}\left(
\begin{array}{l}
k-1 \\
e(\sigma )
\end{array}
\right) } \int_{\Delta^k [0,t]} \circ dB^{\sigma^{-1} \cdot I}.
\]
Due to the assumption that
\[
(X_{t})_{t \geq 0} =^{law} (\Delta_c X_{\frac{t}{c}})_{t \geq 0},
\]
we deduce that
\[
\left( \exp \left( \sum_{k = 1}^{+\infty} \sum_{I=(i_1,...,i_k)}
\Lambda_I (B)_t V_I \right) \right)_{t \ge 0}=^{law} \left(\exp
\left( \sum_{k = 1}^{+ \infty} \sum_{I=(i_1,...,i_k)} \Lambda_I
(B)_{\frac{t}{c}} (\delta_c V)_I \right) \right)_{t \ge 0}.
\]
But since the group $\mathbf{G}$ is nilpotent and simply connected
the exponential map  is a diffeomorphism, therefore
\[
\left(  \sum_{k = 1}^{+\infty} \sum_{I=(i_1,...,i_k)} \Lambda_I
(B)_t V_I  \right)_{t \ge 0}=^{law} \left( \sum_{k = 1}^{+ \infty}
\sum_{I=(i_1,...,i_k)} \Lambda_I (B)_{\frac{t}{c}} (\delta_c V)_I
\right)_{t \ge 0}.
\]
Let us now observe that due to the scaling property of the
fractional Brownian motion
\[
\left( \sum_{k = 1}^{+\infty} \sum_{I=(i_1,...,i_k)} \Lambda_I
(B)_{\frac{t}{c}} (\delta_c V)_I \right)_{t \ge 0} =^{law} \left(
\sum_{k = 1}^{+\infty} \sum_{I=(i_1,...,i_k)} \Lambda_I (B)_{t}
\frac{1}{c^{H\mid I \mid}} \delta_c V_I \right)_{t \ge 0},
\]
where $\mid I \mid$ is the length of the word $I$. Thus, for every
$c>0$,
\[
\left(  \sum_{k = 1}^{+\infty} \sum_{I=(i_1,...,i_k)} \Lambda_I
(B)_t V_I  \right)_{t \ge 0}=^{law}\left( \sum_{k = 1}^{+\infty}
\frac{1}{c^{kH}} \delta_c \left( \sum_{I=(i_1,...,i_k)} \Lambda_I
(B)_{t}  V_I \right) \right)_{t \ge 0}
\]
Let us now observe that $V_1,...,V_d$ is a basis of
$\mathfrak{g}$, therefore all commutators are linear combinations
of the $V_i$'s. By letting $c \rightarrow + \infty$ in the
previous equality, we deduce that $c^{-H} \delta_c$ has to be
bounded when $ c \rightarrow +\infty$. Since $\delta_c =
\text{Exp} ( \phi \ln c)$, all eigenvalues of $\phi$ have a real
part that is smaller than $H$. Finally, by letting $c \rightarrow
0$, we conclude that almost surely, for $k \ge 2$,
\[
\sum_{I=(i_1,...,i_k)} \Lambda_I (B)_{t}  V_I =0.
\]
In particular,
\[
\sum_{1 \le i <j \le d} \left( \int_0^t B^i_s dB^j_s -B^j_s dB^i_s
\right) [V_i,V_j]=0.
\]

Therefore, from the support theorem of \cite{Millet-Sanz} (see
also \cite{coutin-friz-victoir}), all the brackets have to be 0,
that is $\mathbf{G}$ is commutative.
\end{proof}

\begin{remark}
The previous theorem in particular applies to the case of a
Brownian motion, that corresponds to $H=\frac{1}{2}$. In that
case, the theorem can be more easily understood in the following
way (see also \cite{Kunita1} and \cite{Kunita2}). If we by denote
$p_t$ the density of $X_t$ with respect to the Haar measure; it is
well-known that we have an asymptotic development of the form
\[
p_t (1_{\mathbf{G}})\sim_{t \rightarrow 0} \frac{1}{(2\pi
t)^{\frac{d}{2}}}\left(1+ \sum_{k=1}^{+\infty} a_k t^k\right).
\]
By homogeneity, the existence of a scaling property for the
Brownian motion, would then imply $a_k=0$ for every $k$. This
implies that the Riemannian curvature of $\mathbf{G}$ is zero, and
hence $\mathbf{G}$ is commutative.
\end{remark}

If we relax the assumption that the family $(V_1,...,V_d)$ forms a
basis of the Lie algebra $\mathfrak{g}$, we can have a global
scaling property in slightly more general groups than the
commutative ones. Let us first look at one example.

The Heisenberg group $\mathbb{H}$ is the set of $3\times3$
matrices:
\[
\left(
\begin{array}
[c]{ccc}
~1~ & ~x~   & ~z ~\\
~0~ & ~1~   & ~y ~\\
~0~ & ~0~   & ~1 ~
\end{array}
\right)  ,\text{ \ }x,y,z\in\mathbf{R}.
\]
The Lie algebra of $\mathbb{H}$ is generated by the matrices
\[
D_{1}=\left(
\begin{array}
[c]{ccc}
~0~ & ~1~ & ~0~\\
~0~ & ~0~ & ~0~\\
~0~ & ~0~ & ~0~
\end{array}
\right)  ,\text{ }D_{2}=\left(
\begin{array}
[c]{ccc}%
~0~ & ~0~ & ~0~\\
~0~ & ~0~ & ~1~\\
~0~ & ~0~ & ~0~
\end{array}
\right)  \text{ and }D_{3}=\left(
\begin{array}
[c]{ccc}%
~0~ & ~0~ & ~1~\\
~0~ & ~0~ & ~0~\\
~0~ & ~0~ & ~0~
\end{array}
\right)  ,
\]
for which the following equalities hold
\[
\lbrack D_{1},D_{2}]=D_{3},\text{ }[D_{1},D_{3}]=[D_{2},D_{3}]=0.
\]
Consider now the solution of the equation
\[
dX_t = X_t (D_1 dB^1 _t + D_2 dB^2_t),\quad X_0=1,
\]
where $(B^1_t,B^2_t)$ is a two-dimensional fractional Brownian
motion with Hurst parameter $H>\frac{1}{4}$. It is easily seen
that
\[
X_t=\left(
\begin{array}
[c]{ccc}%
~1~ & ~B^1_t~ & ~\frac{1}{2} \left( B^1_t B^2_t +\int_0^t B^1_s dB^2_s -B^2_s dB^1_s\right)~\\
~0~ & ~1~ & ~B^2_t~\\
~0~ & ~0~ & ~1~
\end{array}
\right).
\]
Therefore $(X_{ct})_{t \geq 0} =^{law} (\Delta_c X_{t})_{t \geq
0}$, where $\Delta_c$ is defined by
\[
\Delta_c \left(
\begin{array}
[c]{ccc}
~1~ & ~x~   & ~z ~\\
~0~ & ~1~   & ~y ~\\
~0~ & ~0~   & ~1 ~
\end{array}
\right)= \left(
\begin{array}
[c]{ccc}
~1~ & ~c^H x~   & ~c^{2H} z ~\\
~0~ & ~1~   & ~c^H y ~\\
~0~ & ~0~   & ~1 ~
\end{array}
\right).
\]
In that case, we thus have a global scaling property whereas
$\mathbb{H}$ is of course not commutative but step-two nilpotent.
Actually, we shall have a global scaling property in the Lie
groups that are called the Carnot groups. Let us recall the
definition of a Carnot group.
\begin{definition}
A Carnot group of step (or depth) $N$ is a simply connected Lie
group $\mathbb{G}$ whose Lie algebra can be written
\[
\mathcal{V}_{1}\oplus...\oplus \mathcal{V}_{N},
\]
where
\[
\lbrack \mathcal{V}_{i},\mathcal{V}_{j}]=\mathcal{V}_{i+j}
\]
and
\[
\mathcal{V}_{s}=0,\text{ for }s>N.
\]
\end{definition}

\begin{example}
Consider the set $\mathbb{H}_n =\mathbb{R}^{2n} \times \mathbb{R}$
endowed with the group law
\[
(x,\alpha) \star (y, \beta)=\left( x+y, \alpha + \beta
+\frac{1}{2} \omega (x,y) \right),
\]
where $\omega$ is the standard symplectic form on
$\mathbb{R}^{2n}$, that is
\[
\omega(x,y)= x^t \left(
\begin{array}{ll}
0 & -\mathbf{I}_{n} \\
\mathbf{I}_{n} & ~~~0
\end{array}
\right) y.
\]
On $\mathfrak{h}_n$ the Lie bracket is given by
\[
[ (x,\alpha) , (y, \beta) ]=\left( 0, \omega (x,y) \right),
\]
and it is easily seen that
\[
\mathfrak{h}_n=\mathcal{V}_1 \oplus \mathcal{V}_2,
\]
where $\mathcal{V}_1 =\mathbb{R}^{2n} \times \{ 0 \}$ and
$\mathcal{V}_2= \{ 0 \} \times \mathbb{R}$. Therefore
$\mathbb{H}_n$ is a Carnot group of depth 2 and observe that
$\mathbb{H}_1$ is isomorphic to the Heisenberg group.
\end{example}

Notice that the vector space $\mathcal{V}_{1}$, which is called
the basis of $\mathbb{G}$, Lie generates $\mathfrak{g}$, where
$\mathfrak{g}$ denotes the Lie algebra of $\mathbb{G}$. Since
$\mathbb{G}$ is step $N$ nilpotent and simply connected, the
exponential map is a diffeomorphism. On $\mathfrak{g}$ we can
consider the family of linear operators $\delta_{t}:\mathfrak{g}
\rightarrow \mathfrak{g}$, $t \geq 0$ which act by scalar
multiplication $t^{i}$ on $\mathcal{V}_{i} $. These operators are
Lie algebra automorphisms due to the grading. The maps $\delta_t$
induce Lie group automorphisms $\Delta_t :\mathbb{G} \rightarrow
\mathbb{G}$ which are called the canonical dilations of
$\mathbb{G}$. Let us now take a basis $U_1,...,U_d$ of the vector
space $\mathcal{V}_1$. The vectors $U_i$'s can be seen as left
invariant vector fields on $\mathbb{G}$ so that we can consider
the following stochastic differential equation on $\mathbb{G}$:
\begin{equation}
\label{SDEcarnot} dY_t =\sum_{i=1}^d \int_0^t U_i (Y_s)  dB^i_s,
\text{ } t \geq 0,
\end{equation}
which is easily seen to have a unique solution associated with the
initial condition $Y_0=1_{\mathbb{G}}$. We have then the following
global scaling property:
\begin{proposition}
\[
(Y_{ct})_{t \geq 0} =^{law} (\Delta_{c^H} Y_{t})_{t \geq 0}.
\]
\end{proposition}
\begin{proof}
We keep the notations introduced before the proof of Theorem
\ref{fbmnilpotent}. From Theorem \ref{fbmnilpotent}, we have
\[
Y_{t}= \exp \left( \sum_{k = 1}^{N} \sum_{I=(i_1,...,i_k)}
\Lambda_I (B)_t U_I \right), \text{ } t \geq 0.
\]
Therefore,
\[
\left( Y_{ct} \right)_{t \ge 0} =^{law} \left( \exp \left( \sum_{k
= 1}^{N} c^{H \mid I \mid} \sum_{I=(i_1,...,i_k)} \Lambda_I (B)_t
U_I \right) \right)_{t \ge 0}.
\]
Since
\begin{align*}
\exp \left( \sum_{k = 1}^{N} c^{H \mid I \mid}
\sum_{I=(i_1,...,i_k)} \Lambda_I (B)_t U_I \right)& =\exp \left(
\sum_{k = 1}^{N}  \sum_{I=(i_1,...,i_k)}
\Lambda_I (B)_t (\delta_{c^H} U_I )\right) \\
&=\Delta_{c^H}\exp \left( \sum_{k = 1}^{N}  \sum_{I=(i_1,...,i_k)}
\Lambda_I (B)_t  U_I \right),
\end{align*}
we conclude
\[
(Y_{ct})_{t \geq 0} =^{law} (\Delta_{c^H} Y_{t})_{t \geq 0}.
\]
\end{proof}

The previous proposition admits a counterpart.

\begin{theorem}
Let $\mathbf{G}$ be a connected Lie group (of matrices) with Lie
algebra $\mathfrak{g}$. Let $V_1,...,V_d$ be a family of
$\mathfrak{g}$. Consider now the solution of the equation
\[
dX_t= X_t \left( \sum_{i=1}^d V_i dB^i_t\right), \quad
X_0=1_{\mathbf{G}}.
\]
Assume that there exists a family of Lie group automorphisms
$\Delta_t :\mathbf{G} \rightarrow \mathbf{G}$, $t>0$,  such that:
\begin{enumerate}
\item The map $t \rightarrow \Delta_t$ is continuous and $\lim_{t \rightarrow 0} \Delta_t=1_{\mathbf{G}}$;
\item For $t_1,t_2 \ge 0$, $\Delta_{t_1 t_2} = \Delta_{t_1 } \circ \Delta_{t_2}$;
\item $(X_{ct})_{t \geq 0} =^{law} (\Delta_c X_{t})_{t \geq 0}$;
\end{enumerate}

Then the Lie subgroup $\mathbf{H}$ that is generated by
$e^{V_1},...,e^{V_d}$ is a Carnot group.
\end{theorem}

\begin{proof}

We can follow the lines of the proof of Theorem
\ref{theoremescaling} to deduce that $\mathbf{H}$ has to be a
simply connected nilpotent group and that for every $c>0$,
\[
\left(  \sum_{k = 1}^{+\infty} \sum_{I=(i_1,...,i_k)} \Lambda_I
(B)_t V_I  \right)_{t \ge 0}=^{law}\left( \sum_{k = 1}^{+\infty}
\frac{1}{c^{kH}} \delta_c \left( \sum_{I=(i_1,...,i_k)} \Lambda_I
(B)_{t}  V_I \right) \right)_{t \ge 0},
\]
where $\delta_c$ is the differential map of $\Delta_c$ at
$1_{\mathbf{G}}$. For $k \ge 1$, we denote $\mathcal{V}_k$ the
linear space generated by the set of commutators:
\[
\left\{ V_I, \mid I \mid=k \right\}.
\]
By letting $c \rightarrow +\infty$ and $c \rightarrow 0$ and by
using the support theorem of \cite{coutin-friz-victoir}, we obtain
that $c^{-kH} \delta_c$ is bounded on $\mathcal{V}_k$ for $c
\rightarrow +\infty$ and $c \rightarrow 0$. Since $c^{-kH}
\delta_c = \text{Exp} ( (\phi-kH \mathbf{Id}) \ln c)$, for some
matrix $\phi$, we conclude
\[
\mathfrak{h}=\bigoplus_{k=1}^{+\infty} \mathcal{V}_k,
\]
where $\mathfrak{h}$ is the Lie algebra of $\mathbf{H}$. This
proves that $\mathbf{H}$ is a Carnot group.
\end{proof}

\section{Stochastic analysis of the fractional Brownian motion on a Lie group}

As before, let $\mathbf{G}$ be a compact connected Lie group of
matrices with Lie algebra $\mathfrak{g}$. From the theory of
compact Lie groups, it is well-known that $\mathbf{G}$ can be
endowed with a bi-invariant Riemannian structure. This
bi-invariant Riemannian structure induces a scalar product
$\langle \cdot , \cdot \rangle_{\mathfrak{g}}$ on $\mathfrak{g}$
for which the maps
\begin{align*}
\mathbf{Ad}_g : & \mathfrak{g} \rightarrow \mathfrak{g} \\
                & x  \rightarrow  g x g^{-1},
\end{align*}
are isometries. The map $ g \rightarrow \mathbf{Ad}_g$ is called
the adjoint representation of $\mathbf{G}$ on $\mathfrak{g}$.

Let $V_1,...,V_d$ be a family of $\mathfrak{g}$ that forms an
orthonormal basis for the scalar product $\langle \cdot , \cdot
\rangle_{\mathfrak{g}}$.

Consider now the solution of the equation
\[
dX_t= X_t \left( \sum_{i=1}^d V_i dB^i_t\right), \quad
X_0=1_{\mathbf{G}}.
\]
where $(B_t)_{t \ge 0}$ is a $d$-dimensional fractional Brownian
motion with Hurst parameter $H > \frac{1}{2}$. Our goal in this
section is to show that for every $t >0$ the random variable $X_t$
admits a density with respect to the Haar measure of $\mathbf{G}$.
For that, we shall develop a
  Malliavin calculus for the fractional Brownian motion on the path group
space (see \cite{Malliavin} for the Brownian case).

Let us recall (see for instance \cite{du}) that the fractional
Brownian motion $(B_t)_{t \ge 0}$ admits the Volterra type
representation:
\[
B^i_t=\int_0^t \mathbf{K}(t,s)dW^i_s,
\]
where $(W_t)_{t \ge 0}$ is a $d$-dimensional Brownian motion and
\begin{align*}
\mathbf{K} (t, s)=c_H s^{\frac{1}{2} - H} \int_s^t
(u-s)^{H-\frac{3}{2}} u^{H-\frac{1}{2}} du, \text{ }t>s,
\end{align*}
where $c_H =\sqrt{ \frac{ H(2H-1)}{\beta (2-2H,H-\frac{1}{2})}} $
is the constant such that the covariance function is
\[
R\left( t,s\right) =\frac{1}{2}\left(
s^{2H}+t^{2H}-|t-s|^{2H}\right).
\]
The key-point of our study is the following integration by parts
formula.

\begin{proposition}
Let $f : \mathbf{G} \rightarrow \mathbb{R}$ and $h :
\mathbb{R}_{\ge 0} \rightarrow \mathfrak{g}$ be smooth functions,
then we have
\[
\mathbb{E} \left( \left\langle \nabla f (X_t) , \int_0^{t} \left(
\int_s^t \frac{\partial \mathbf{K}}{\partial u} (u,s)
\mathbf{Ad}_{X_u} du \right) h'(s)ds \right\rangle_{\mathfrak{g}}
\right)=\mathbb{E} \left( f (X_t) \int_0^t \left\langle h'(s)
,\sum_{i=1}^d V_i dW^i_t \right \rangle_{\mathfrak{g}} \right),
\]
where $\nabla f= \sum_{i=1}^d (V_i f) V_i$.
\end{proposition}

\begin{proof}
Let $\varepsilon >0$. We denote
\begin{equation*}
W^{\varepsilon}=W+\varepsilon \int_0^{\cdot} h'(s) ds.
\end{equation*}
With obvious notations, we have then
\begin{equation*}
B^{\varepsilon}=B + \varepsilon \int_0^{\cdot} \mathbf{K}
(\cdot,s)h'(s) ds.
\end{equation*}
If we consider now the stochastic development on $\mathbf{G}$ of
$B^{\varepsilon}$, that is the solution of the equation
\[
dX^{\varepsilon}_t= X^{\varepsilon}_t \left( \sum_{i=1}^d V_i
dB^{i,\varepsilon}_t\right), \quad
X^{\varepsilon}_0=1_{\mathbf{G}},
\]
we claim that for any smooth $f : \mathbf{G} \rightarrow
\mathbb{R}$,
\begin{equation}\label{variationinfinitesimale}
\left( \frac{d}{d\varepsilon} \right)_{\varepsilon = 0} f (
X_t^{\varepsilon}) =\left\langle \nabla f (X_t) , \int_0^{t}
\left( \int_s^t \mathbf{Ad}_{X_u} \frac{\partial
\mathbf{K}}{\partial u} (u,s) du \right) h'(s)ds
\right\rangle_{\mathfrak{g}}, \text{ } t>0.
\end{equation}
Indeed let us recall a  well-known fact in the theory of linear
equations:

Let $(x_t)_{t \ge 0}$ and $(\rho_t)_{t \ge 0}$ be piecewise $C^1$
and $\mathbb{R}^d$-valued paths. Denote $(y_t)_{t \ge 0}$ and
$(y^{\varepsilon}_t)_{t \ge 0}$ the respective solutions of the
linear equations:
\[
dy_t= y_t \left( \sum_{i=1}^d V_i dx^i_t\right), \quad
y_0=1_{\mathbf{G}},
\]
and
\[
dy^{\varepsilon}_t= y^{\varepsilon}_t \left( \sum_{i=1}^d V_i
d(x^i_t+\varepsilon \rho^i_t) \right), \quad
y^{\varepsilon}_0=1_{\mathbf{G}}.
\]
Then we have
\[
\left( \frac{dy^{\varepsilon}_t}{d\varepsilon}
\right)_{\varepsilon = 0}=\left( \int_0^t \mathbf{Ad}_{y_s}
\left(\sum_{i=1}^d V_i d\rho^i_s \right) \right)y_t.
\]
By applying this result to piecewise linear interpolations of $B$
along the dyadic subdivision and by passing to the limit with the
transfer principle of \cite{Li-Lyons} we get formula
(\ref{variationinfinitesimale})

We deduce therefore our integration by parts formula from the
classical integration by parts formula on the Wiener space.

\end{proof}

\begin{remark}
As an easy consequence, we can derive a Clark-Ocone type formula
for $(X_t)_{t \ge 0}$: If $f : \mathbf{G} \rightarrow \mathbb{R}$
is a smooth function, then we have
\[
f(X_t) = \mathbb{E} \left( f(X_t) \right) + \int_0^t \int_s^t
\left\langle \frac{\partial \mathbf{K}}{\partial u} (u,s)
\mathbb{E} \left( (\mathbf{Ad}_{X_u} )^{\ast} \nabla f (X_t) \mid
\mathcal{F}_s \right) du, \sum_{i=1}^d V_i dW^i_s \right
\rangle_{\mathfrak{g}}.
\]
\end{remark}

\begin{remark}
More generally, let $( \mathbb{M},g)$ be a $d$-dimensional compact
connected Riemannian manifold. We denote $\mathcal{O}\left(
\mathbb{M}\right) $ the orthonormal frame bundle of $\mathbb{M}$
and $ \pi $ the canonical surjection on $\mathbb{M}$. Let $\left(
H_{i}\right) _{i=1,...,d}$ the horizontal fundamental vector
fields of $\mathcal{O}\left( \mathbb{M}\right) $ and
$\overline{\Omega}$ the equivariant representation of the Cartan
curvature form $\Omega$. Consider the rough paths differential
equation
\[
dX^*_t=  \sum_{i=1}^d H_i (X^*_t ) dB^i_t, \quad X^*_0=m^* \in
\mathcal{O}\left( \mathbb{M}\right),
\]
Due to the compactness of $\mathbb{M}$, from rough paths theory
\cite{LyQi} and \cite{CQ}, it is easy to show that the previous
equation has a unique solution. By following the
Eels-Elworthy-Malliavin's approach to  Riemannian Brownian motion
(see \cite{Malliavin}), one can define a fractional Brownian
motion $(X_t)_{t \ge 0}$ on $\mathbb{M}$ by:
\[
(X_t)_{t \ge 0}=(\pi X_t^*)_{t \ge 0}.
\]
In this setting, by using Malliavin-Cruzeiro isomorphism (see
\cite{Malliavin}), we can show the following integration by parts
formula. Let $h=(h_1,...,h_d): \mathbb{R}_{\ge 0} \rightarrow
\mathbb{R}^d$ be a smooth function. For any smooth $\xi :
\mathbb{M} \rightarrow \mathbb{R}$,
\begin{equation*}
\mathbb{E} \left( \left\langle \nabla \xi \left( X_{t} \right) ,
X_t^* \Theta_t \right\rangle_{T_{X_t}
\mathbb{M}}\right)=)=\mathbb{E} \left( \xi (X_t)\sum_{i=1}^d
\int_0^t h'_i (s)dW^i_s \right),  \text{ }t>0,
\end{equation*}
where $\Theta$ is the $\mathbb{R}^d$-valued process  solution of
the rough paths differential equation
\begin{equation*}
\Theta_t =  \int_0^t \mathbf{K} (t,s) h'(s) ds+\int_0^t  \int_0^s
\overline{\Omega}_{X_u^*} ( \Theta_u , dB_u).dB_s.
\end{equation*}
\end{remark}

As a corollary of the previous integration by parts formula, we
deduce:

\begin{proposition}
For every $t >0$ the random variable $X_t$ admits a density with
respect to the Haar measure of $\mathbf{G}$.
\end{proposition}

\begin{proof}
We have to show that for $t >0$, the Malliavin matrix $\Gamma_t$
is invertible. Here, we have from the integration by parts
formula:
\[
\Gamma_t=\int_0^t \left( \int_s^t \frac{\partial
\mathbf{K}}{\partial u} (u,s) \mathbf{Ad}_{X_u} du \right)^*
\left( \int_s^t \frac{\partial \mathbf{K}}{\partial u} (u,s)
\mathbf{Ad}_{X_u} du \right) ds
\]
If $x \in \mathbf{Ker} \text{ }\Gamma_t$, then we have
\[
x^* \Gamma_t x=0.
\]
It implies that for $s \le t$,
\[
\int_s^t \frac{\partial \mathbf{K}}{\partial u} (u,s)
\mathbf{Ad}_{X_u} (x) du=0,
\]
and so for $s \le t$,
\[
\mathbf{Ad}_{X_s} (x)=0.
\]
Since $\mathbf{Ad}_{X_s}$ is an isometry, we get therefore $x=0$.
\end{proof}

\end{document}